\documentclass{article}

\usepackage[latin1]{inputenc}
\usepackage{graphicx}
\usepackage{epsfig}

\newenvironment{proof}{\trivlist
\item[\hskip\labelsep{\it Proof}\,:]}{\hfill{$q.e.d$}\endtrivlist}
\newtheorem{theorem}{Theorem}
\newtheorem{definition}{Definition}

\newtheorem{lemma}{Lemma}
\newtheorem{remark}{Remark}
\newtheorem{corollary}{Corollary}

\newfont{\bb}{msbm10 at 11pt}
\def\r{\hbox{\bb R}}

\begin{document}

\title{Stationary  rotating surfaces in Euclidean space}
\author{Rafael L\'opez\footnote{Partially
supported by MEC-FEDER
 grant no. MTM2007-61775.
}\\
 Departamento de Geometr\'{\i}a y Topolog\'{\i}a\\
Universidad de Granada\\
18071 Granada (Spain)\\
e-mail: {\tt rcamino@ugr.es}}

\date{}

\maketitle

\begin{abstract} A stationary rotating surface is a   compact surface in Euclidean space whose mean curvature $H$ at each point $x$ satisfies $2H(x)=a r^2+b$, where $r$ is the distance from $x$ to a fixed straight-line $L$, and $a$ and $b$ are constants. These surfaces are solutions of a variational problem that describes the shape of a drop of incompressible fluid in equilibrium by the action of surface tension when it rotates about $L$ with constant angular velocity. The effect of gravity is neglected. In this paper we study the geometric configurations of such surfaces, focusing the relationship  between the geometry of the surface and the one of its boundary. As special cases, we will consider two families of such surfaces: axisymmetric surfaces and embedded surfaces with planar boundary.

\end{abstract}


\section{Introduction}


In the absence of gravity, we consider the steady rigid rotation of
an homogeneous incompressible fluid drop which is surrounded by a
rigidly rotating incompressible fluid. Our interest is the study of
the shape of such drop when it attains a state of mechanical
equilibrium. In such case, we will call it a
rotating liquid drop, or simply, a rotating drop. Rotating liquid drops have been the
subject of intense study beginning from the work of Joseph Plateau
\cite{pl}. Experimentally, he observed  a variety of axisymmetric
shapes that can summarized as follows: starting with zero angular
velocity, we begin with a spherical shape. As we increase the angular velocity,
the drop changes through a sequence of shapes which evolved from
axisymmetric for slow rotation to ellipsoidal and two-lobed and
finally toroidal at very large rotation. This was theoretically
shown by Poincar\'e \cite{po}, Beer \cite{be},  Chandrasekhar \cite{ch} and Brown and
Scriven \cite{bs1,bs2}.  The experiments of Plateau inspired  an scientific  interest since
they could be models in other areas of physics, such as,
astrophysics, nuclear physics, fluid dynamics, amongst others. For example,  they arise in celestial mechanics in the study of self-rotating stars and planets \cite{ch2,ko,wav}. Experimentally, these figures have appeared in microgravity environments and experimental works in absence of gravity
\cite{la,wa}. From a theoretical viewpoint,  there is a mathematical interest for rotating drops, focusing in subjects such as existence and stability. The literature  is extensive, and without trying to give a complete list, we refer the reader to \cite{ag,at,au,cf,co,ceg,sr}

Let us take usual coordinates $(x_1,x_2,x_3)$ and assume the liquid drop rotates about the
$x_3$-axis with a  constant angular velocity $\omega$. Let $\rho$ be the constant density of the fluid of the drop.
Let $W$ be the bounded open set
in $\r^3$, which is the region occupied by the rotating drop. We set $S=\partial W$ as the free interface between the drop and the ambient liquid  and that we suppose to be a smooth boundary surface. The energy of this mechanical system is given by
$$E=\tau|S|-\frac12\rho\omega^2\int_{W} r^2 dx,$$
where $\tau$ stands for the surface tension on $S$, $|S|$ is the surface area of $S$ and $r=r(x)=\sqrt{x_1^2+x_2^2}$ is the distance
from a point $x$ to the $x_3$-axis, the axis of rotation.
The term  $\tau|S|$ is the surface energy of the drop and
$\frac12\rho\omega^2\int_W r^2 dx$ is the potential energy associated
with the centrifugal force. We assume that the volume $V$ of the drop remains constant while  rotates.

We seek the shape of the liquid drop when the
configuration is stationary, that is, the drop is a critical point of the energy for all volume preserving perturbations.
The equilibrium is obtained by the balance between the capillary
force that comes from the surface tension of $S$  and the centrifugal force of the rotating liquid.
The equilibrium shapes of such a drop are governed by the
Young-Laplace equation
$$2\tau H(x)=-\frac12\rho\omega^2 r^2+\lambda\hspace*{1cm}\mbox{$x\in S$},$$
where $\lambda$ is a constant depending on the volume constraint.
As consequence, the mean curvature of the interface $S$ satisfies an equation of type
$$H(x_1,x_2,x_3)=a r^2+b,$$
where $a,b\in\r$. We say then that $S$ is a stationary rotating surface.

In the case that the interface $S$ is an embedded surface (no self-intersections), Wente showed that a rotating drop has a plane of symmetry perpendicular to
the  $x_3$-axis and any line parallel to the axis and meeting the drop cuts it in a segment whose center lies on the plane of symmetry \cite{we2}. Moreover, the plane  through the mass centre perpendicular to the axis of the rotation coincides with the plane of symmetry of the rotating liquid drop.

The first configurations studied in the literature are the axisymmetric shapes of rotating liquid drops, that is,  surfaces of revolution with respect to the axis of rotation.  In such case, the Young-Laplace equation is a second order differential equation and a first integration can done (see Section \ref{sect-axi}). The purpose of this paper is to present the study of rotating liquid drops in a more general sense, assuming for example that the surface is not rotational, nor embedded or with a possible non-empty boundary.

When the liquid does not rotate, that is, the angular velocity is zero, $\omega=0$ (or $a=0$ in the Laplace equation),   the mean curvature of the surface is constant and we abbreviate  by saying a cmc-surface. Although in this paper we discard this situation, stationary rotating surfaces share techniques and type of results with cmc-surfaces.  Actually
part of our work follows the same scheme and methodology, although the Laplace equation in our setting is more difficult and the results are less definitive. This can clearly see in the case that the surface is embedded. In this sense, it is worthwhile saying two facts that  makes  different both settings:
\begin{enumerate}
\item There are stationary rotating surfaces with genus $0$ that are not embedded. See Figure \ref{fig2} (right). However, the celebrated Hopf's theorem \cite{ho} asserts that round spheres are the only cmc-surfaces with genus $0$.
\item There are toroidal rotating drops that are embedded. In the family of  cmc embedded closed surfaces, the only possibility is the round sphere (Alexandrov's theorem \cite{al}).
\end{enumerate}
Since it is rather difficult to consider the case in which the liquid does not form a surface of revolution  and because axisymmetric shapes are more suitable to study, a first question is whether a stationary rotating embedded closed surface must be a surface of revolution.  As we have mentioned, Wente's theorem assures that there exists a plane of symmetry perpendicular to the axis of rotation. For this, he  utilises the so-called Alexandrov reflection method  by horizontal planes. However, one cannot do a similar argument with vertical planes: in general, when one reflects the surface about these planes,  it cannot compare the value of the mean curvature at the contact points, on the contrary what occurs for cmc-surfaces and capillary surfaces \cite{al,we3}.

Assume now that the boundary of the surface is a non-empty set. The simplest case to consider is that the boundary is a horizontal circle centred at the origin. A natural question is  whether
a compact rotating liquid drop in $\r^3$  bounded by a circle
is necessarily a rotational surface. More generally, one can consider the problem whether a stationary rotating surface inherits the symmetries of its boundary.

In Section \ref{sect-preli}, we formulate the Laplace equation and we derive the second variation of the energy. In Section  \ref{sect-integral}, a set of integral formulae will be obtained relating quantities between the surface and its boundary. Section \ref{sect-axi} considers axisymmetric rotating surfaces obtaining some estimates of the profile curve, the area and the  volume of the surface. Finally, sections \ref{sect-embedded} and \ref{sect-estable} are devoted to analyse configurations and stability of liquid drops with planar boundary.

\section{Preliminaries}\label{sect-preli}

Let $\r^3$ be the Euclidean three-space and let
$(x_1,x_2,x_3)$ be the usual coordinates. Let $M$  be an oriented (connected) compact surface  and we shall
denote by $\partial M$ the boundary of $M$.
Consider a smooth immersion $x:M\rightarrow\r^3$ and let  $N$ be the  Gauss map. Denote by $\{E_1,E_2,E_3\}$ the canonical orthonormal base in $\r^3$:
$$E_1=(1,0,0),\hspace*{.5cm}E_2=(0,1,0),\hspace*{.5cm}E_3=(0,0,1).$$
We write $N_i=\langle N,E_i\rangle$ and $x_i=\langle x,E_i\rangle$
for $1\leq i\leq 3$.  We define the energy functional as
$$E(x)=\tau\int_M dM-\frac12\rho\omega^2\int_M  r^2 x_3 N_3\ dM,\hspace*{1cm}r=\sqrt{x_1^2+x_2^2}$$
where $\tau,\rho$ and $\omega$ are constant and $dM$ is the area element on $M$. Here $\int_M  r^2 x_3 N_3\ dM$ represents the centrifugal force of  the surface with respect to the $x_3$-axis.
Consider a smooth variation $X$ of $x$, that is, a smooth map $X:(-\epsilon,\epsilon)\times M:M\rightarrow\r^3$  such that, by setting $x_t=X(t,-)$, we have $x_0=x$ and $x_t-x\in C^{\infty}_0(M)$. Let $u\in C^{\infty}(M)$ be the normal component of the variational vector field of $x_t$,
$$u=\langle\frac{\partial x_t}{\partial t}{\bigg|}_{t=0},N\rangle.$$
We take $E(t):=E(x_t)$ the value of the energy for each immersion $x_t$. The first variation of $E$ at $t=0$ is given by
$$E'(0)=-\int_M\bigg(\frac12\rho\omega^2 r^2+2\tau H\bigg) u\ dM.$$
Here $H$ is the mean curvature of the immersion $x$.
We require that the volume $V(t)$ of each immersion $x_t$  remains constant throughout the variation.  The first variation of the volume functional is
$$V'(0)=\int_M u\ dM.$$
By the method of Lagrange  multipliers,  the first variation of $E$ at $t=0$ is to
be zero relative for all volume preserving variations if there is a constant $\lambda$ so that
$E'(0)+\lambda V'(0)=0$.  This yields the
condition
$$2\tau H=-\frac12\rho\omega^2 r^2+\lambda\hspace*{1cm}\mbox{on $M$}.$$
See \cite{we2} for details.  Thus the mean curvature $H$ satisfies an equation of type
\begin{equation}\label{laplace}
2H(x)=a r^2+b,\hspace*{1cm}a,b\in\r.
\end{equation}

\begin{definition} Let $L$ be a straight-line of $\r^3$. A stationary rotating surface (with respect to $L$) is an oriented compact surface $M$ immersed in $\r^3$ such that the mean curvature of the immersion satisfies the Laplace equation (\ref{laplace}), where $r=\mbox{dist}(x,L)$. If the surface is embedded, we say that $M$ is a rotating liquid drop.
\end{definition}

\begin{remark} Throughout this work, we suppose that the straight-line $L$  is the $x_3$-axis. Thus, $r^2=x_1^2+x_2^2$. Moreover, we shall use the words "horizontal" and "vertical" with respect to $L$, that is, by "horizontal" we mean orthogonal to $L$ and by "vertical", we mean parallel to $L$.
\end{remark}

We need to precise the definition that a fixed curve of $\r^3$ is the boundary of an immersion $x:M\rightarrow\r^3$. Let $\Gamma$ be a closed curve in $\r^3$. We say the $\Gamma$ is the boundary of $M$ if  $x_{|\partial M}:\partial M\rightarrow\r^3$ is an embedding, with $x(\partial M)=\Gamma$.

Given the definition of a stationary rotating surface, we may derive  the second variation of the energy of critical points in order to give the notion of stability. Stability of rotating liquid drops has been studied for rotationally symmetric configurations in \cite{agu,ca,mc}. A general formula of the second variation was obtained by Wente \cite{we} (see also \cite{mc}). We give a different method for this computation following  ideas of M. Koiso and B. Palmer \cite{kp}.  Assume that $x$ is a critical point of $E$ and we calculate the second variation of the functional $E$. For this, we write $E''(0)$ in the form
$$E''(0)=-\int_M u\cdot L[u]\ dM,$$
where $L$ is a linear differential operator acting on  the normal component, which we want to find it.
The computation of the operator $L$ is as follows.
Since the translations in the $E_3$ are symmetries of the energy functional $E$, then $L[N_3]=0$. Moreover,  the rotation with respect to the $x_3$-axis is a  symmetry of the functional and then $L[\psi]=0$, where the function $\psi$ is $\psi=\langle x\wedge N,E_3\rangle$ and $\wedge$ is the vector product of $\r^3$. We now compute  $L[N_i]$ and $L[\psi]$. The tension field of the Gauss map satisfies
$$\Delta N+|\sigma|^2 N=-2\nabla H,$$
where $\Delta$ is the Laplacian in the metric induced by $x$,
$\sigma$ is the second fundamental of the immersion  and $\nabla$ is
the covariant differentiation. Since $2H=ar^2+b$, we have
\begin{equation}\label{dh}
2\nabla H=2a(x_1\nabla x_1+x_2\nabla x_2)=2a\bigg( x_1 E_1+x_2 E_2-(h-x_3 N_3)N\bigg).
\end{equation}
Here $h=\langle N,x\rangle$ stands for the support function of $M$. Thus
\begin{equation}\label{ene}
\Delta N+\bigg(|\sigma|^2-2a(h-x_3 N_3)\bigg)N=-2a (x_1 E_1+x_2 E_2).
\end{equation}
In particular,
\begin{equation}\label{ene2}
\Delta N_3+\bigg(|\sigma|^2-2a(h-x_3 N_3)\bigg)N_3=0
\end{equation}
We now take the function $\psi$. In general, we have
$$\Delta \psi+|\sigma|^2 \psi=-2\langle \nabla H,E_3\wedge x\rangle.$$
It follows from (\ref{dh}) that
$$-2\langle\nabla H,E_3\wedge x\rangle=2a(h-x_3 N_3)\langle N,E_3\wedge x\rangle=2a(h-x_3 N_3)\psi.$$
Therefore,
\begin{equation}\label{ene3}
\Delta \psi+\bigg(|\sigma|^2-2a(h-x_3 N_3)\bigg)\psi=0.
\end{equation}
As a consequence of (\ref{ene2}) and (\ref{ene3}),
$$L=\Delta+|\sigma|^2-2a(h-x_3 N_3).$$
\begin{definition} Consider $x:M\rightarrow\r^3$ be a smooth immersion that satisfies the Laplace equation (\ref{laplace}).
We say that $x$ is stable if
\begin{equation}\label{second}
-\int_M u\bigg(\Delta u+\big(|\sigma|^2-2a(h-x_3 N_3)\big)u\bigg)\ dM\geq 0
\end{equation}
for all $u\in C_0^{\infty}(M)$ such that
$$\int_M u\ dM=0.$$
The immersion is called strongly stable is (\ref{second}) holds for all $u\in C_0^{\infty}(M)$.
\end{definition}

\section{Integral formulae for stationary rotating surfaces with boundary}\label{sect-integral}

In this section we develop a series of integral formulae for stationary rotating surfaces that will be  used in further sections. As an application, we obtain here an estimate of the height of such surfaces.  We define  the vector valued 1-form $\omega_p= x(p)\wedge v$, $p\in M$, $v\in T_p M$. Then $d\omega=2N$, and the Stokes formulas gives
\begin{equation}\label{n1}
2\int_M N\ dM=\int_{\partial M} x\wedge\alpha'\ ds,
\end{equation}
where $\alpha$ is a parametrization by the length-arc of $\partial M$ that orients $\partial M$ by the induced orientation from $M$ and $ds$ is the length-arc element. Now, consider $\mu_i=x_i^2 \omega$, for $1\leq i\leq 3$. A straightforward computation leads to
$$d\mu_i=4x_i^2 N-2 h x_i E_i.$$
By integrating on $M$, we have
\begin{equation}\label{eq1}
4\int_M x_i^2 N\ dM-2\int_M hx_i E_i\ dM=\int_{\partial M}x_i^2\ x\wedge\alpha'\ ds.
\end{equation}
We now define the  $1$-form  $\beta_p(v)= N(x(p))\wedge v$. Then $d\beta=-2HN$ and Stokes's formula yields
\begin{equation}\label{eq2}
a\int_M r^2 N\ dM+b\int_M N\ dM=-\int_{\partial M}\nu\ ds,
\end{equation}
where $\nu$  the inward   conormal unit vector  field along  $\partial M$. Equation (\ref{eq2}) can also be  obtained by considering the equation
\begin{equation}\label{delta}
\Delta x=2HN=(ar^2+b)N,
\end{equation}
 which holds for any immersion $x$. Then we apply the divergence theorem obtaining (\ref{eq2}) again.

A first consequence of these formulas is the following result about the mass center of a rotating liquid drop (see also \cite{sr}).

\begin{theorem} \label{t-mass} Assume that $M$ is a  stationary rotating embedded closed surface.  Then the mass center of the surface lies at the $x_3$-axis.
\end{theorem}

\begin{proof}
Let $W$ denote  the enclosed domain  by $x(M)$. If $\Delta_0$ denotes the Euclidean Laplacian operator, $\Delta_0 x_i^3=6x_i$, $1\leq i\leq 3$. The divergence theorem gives
\begin{equation}\label{eq3}
2\int_W x_i\ dV=\int_M x_i^2 N_i\   dM,\hspace*{1cm}1\leq i\leq 3,
\end{equation}
where $dV$ is the volume element of $\r^3$. We multiply by $E_j$ in equations (\ref{eq1}) and (\ref{eq2}). Then for  $i,j\in\{1,2\}$ we have
$$2\int_M x_i^2 N_j\ dM=\delta_{ij}\int_M h x_i\ dM,\hspace*{.5cm}\mbox{and}\hspace*{.5cm}\int_M r^2 N_i\ dM=0.$$
As a consequence,
$$\int_M x_i^2 N_j\ dM=0$$
for all $i,j\in\{1,2\}$. Using the above equation together (\ref{eq3}), we obtain
$$\int_W x_i\ dV=0,\hspace*{1cm}i=1,2,$$
which proves the result.

\end{proof}

We extend the above result for rotating liquid drops orthogonally deposited on a horizontal plane.

\begin{theorem}  Let $M$ be a stationary rotating embedded surface. Assume that $\partial M$ is  contained in a horizontal plane $P$. If $M$ lies in one side of $P$ and $M$  is orthogonal to $P$ along $\partial M$, then  the mass center of the surface lies at the $x_3$-axis
\end{theorem}

\begin{proof}
As $\partial M$ is a planar curve and $M$ is orthogonal to $P$, then  $\langle x\wedge\alpha',E_i\rangle=0$ and $\langle\nu, E_i\rangle=0$, respectively,  for $i=1,2$. Multiplying by $E_i$ in (\ref{eq1}) and (\ref{eq2}), we obtain for each $i\in\{1,2\}$
\begin{eqnarray}
4\int_M x_i^2  N_i\ dM-2\int_{M} h x_i\ dM&=& 0, \label{r1}\\
2\int_M x_i^2 N_j\ dM &=&\delta_{ij}\int_M h x_i\ dM \label{r2}\\
a\int_M r^2 N_i\ dM&=&-b\int_M N_i \ dM-\int_{\partial M}\nu_i=0,
\end{eqnarray}
As a consequence of the above three equations, we have
\begin{equation}\label{r4}
\int_M x_i^2 N_i\ dM=0,\hspace*{1cm}i=1,2.
\end{equation}
The calculation of the mass center of the surface follows the same steps than in Theorem \ref{t-mass}. One begins by considering the closed surface $M\cup\Omega$, where $\Omega\subset P$ is the bounded domain by $\partial M$ and let $W$ be the bounded domain of $\r^3$ that determines. If $\eta_\Omega$ is  the induced orientation on $\Omega$, we use (\ref{r4}) and the fact that $\langle\eta_\Omega,E_i\rangle=0$, $i=1,2$, to conclude
$$2\int_W x_i\ dV=\int_M x_i^2 N_i\ dM +\int_\Omega x_i^2\langle \eta_\Omega,E_i\rangle\ d\Omega=\int_M x_i^2 N_i\ dM=0,\hspace*{1cm}i=1,2.$$
The result now follows.
\end{proof}

In the next theorem we obtain an integral formula where all integrals are  evaluated on $\partial M$.

\begin{theorem} Let $M$ be a stationary rotating surface with non-empty boundary. Then
\begin{equation}\label{ar2}
\int_{\partial M}(ar^2+2b)\langle x\wedge\alpha',E_3\rangle\ ds=-4\int_{\partial M}\nu_3\ ds.
\end{equation}
In particular, we have
\begin{equation}\label{ar22}
a\int_{\partial M}r^2\ ds\leq 4\bigg(L(\partial M)-b\ a(\partial M)\bigg),
\end{equation}
where $L(\partial M)$ is the length of $\partial M$ and $a(\partial M)$ is the algebraic area of $\partial M$.
\end{theorem}

\begin{proof} Multiplying by $E_3$ in  (\ref{n1}) and (\ref{eq1}), we know that
\begin{eqnarray}\label{n3}
2 \int_M N_3\ dM&=&\int_{\partial M}\langle x\wedge\alpha',E_3\rangle\ ds.\\
4\int_M x_i^2 N_3\ dM&=&\int_{\partial M} x_i^2 \langle x\wedge\alpha',E_3\rangle\ ds\hspace*{1cm} i=1,2.
\end{eqnarray}
Then
\begin{equation}\label{n32}
4\int_M r^2 N_3\ dM=\int_{\partial M}r^2 \langle x\wedge\alpha',E_3\rangle\ ds.
\end{equation}
From (\ref{eq2}) and (\ref{n3}), we obtain
$$a\int_M r^2 N_3\ dM+\frac{b}{2}\int_{\partial M} \langle x\wedge\alpha',E_3\rangle\ ds=-\int_{\partial M}\nu_3 ds.$$
By combination this equation with (\ref{n32}), we conclude (\ref{ar2}), and this completes the proof.
\end{proof}

Equation (\ref{ar2}) (or (\ref{ar22})) can be viewed as a necessary condition for the existence of a stationary rotating surface with a  prescribed curve as its boundary. Exactly, we propose the following
 \begin{quote} \emph{Problem.} Let $\Gamma$ be a closed  curve in Euclidean space and $a,b\in\r$. Does exist a stationary rotating surface $M$ bounded by $\Gamma$ and with mean curvature $2H(x)=ar^2+b$?
\end{quote}
In general, the answer is "No" because from (\ref{ar2})  it is necessary a certain relation between the quantities $a,b, L(\Gamma)$ and $a(\Gamma)$. For example, in the simplest case of $\Gamma$, that is, a horizontal circle of radius $R$, we have

\begin{corollary} Let $\Gamma$ be a horizontal circle of radius $R$ centred at the $x_3$-axis. If  $M$ is a stationary rotating surface bounded by $\Gamma$ then
\begin{equation}\label{heinz}
|aR^2+2b|\leq\frac{4}{R}.
\end{equation}
\end{corollary}

\begin{proof}
It is sufficient to do $\langle\nu,E_3\rangle\leq 1$ in Equation (\ref{ar2}).
\end{proof}

This corollary has the same flavour than a classical result due to Heinz \cite{he}, which asserts that a necessary condition for the existence of a compact surface with constant mean curvature $H$  bounded by a circle of radius $R$ is that $|H|\leq 1/R$. Actually, and in our setting, if $H$ is constant, then $a=0$ and $2H=b$. Then (\ref{heinz}) reads as $|2b|=|4H|\leq 4/R$, rediscovering the Heinz's result.

This section finishes with an application of formula (\ref{ar2}) in order to derive a height estimate for a stationary rotating graph.

\begin{theorem}\label{height}
Let $M$ be a stationary rotating surface that is a graph on a horizontal plane $P$ and $\partial M\subset P$. Denote $R=\max_{x\in\partial M}r(x)$. Assume that the mean curvature is $2H(x)=ar^2+b$, where $a\not=0$ and $ab\geq 0$. If
$h=\max_{x\in M}\mbox{\rm dist}(x,P)$,
then
\begin{equation}\label{estimate}
h\leq \frac{|a R^2+2b|}{8\pi}\ \mbox{\rm area}(M).
\end{equation}
\end{theorem}

\begin{proof} After a vertical displacement, we assume that $P$ is the plane $x_3=0$ and that $M=\mbox{graph}(u)$, where $u$ is a smooth function on a domain $\Omega\subset P$. Because $a\not=0$ and $ab\geq 0$, the mean curvature  satisfies $H\geq 0$ on $M$ or $H\leq 0$ on $M$. The maximum principle implies that $M$ lies in one side of $P$. If it is necessary, after a reflection about $P$, we can suppose that $M$ lies in the upper half-space determined by $P$ and that the orientation points downwards (this implies $H\geq 0$ and $a,b\geq 0$). Let us introduce the following notation:
$$P_t=\{x\in\r^3;x_3=t\}\hspace*{1cm}M(t)=\{x\in M;x_3\geq t\},\hspace*{1cm}\Gamma(t)=M(t)\cap P_t,$$
and $A(t)$ and $L(t)$ the area and length of $M(t)$ and $\Gamma(t)$ respectively. Let $\Omega(t)$ be the planar domain of $P_t$ bounded by $\Gamma(t)$. We apply Equation (\ref{ar2}) for each surface $M(t)$ obtaining
\begin{equation}\label{coarea}
\int_{\Gamma(t)}\nu^t_3\ ds_t\leq \frac14 \int_{\Gamma(t)}(a r^2+2b) \langle x\wedge\alpha',E_3\rangle\ ds_t\leq \frac{1}{2}|\Omega(t)| (a R^2+2b),
\end{equation}
where $\nu_t$ is the inner conormal unit vector to $M(t)$ along $\Gamma(t)$, $ds_t$ is the induced length arc of $\Gamma(t)$ and $|\Omega(t)|$ is the area of  $\Omega(t)$.
We utilise   the H\"{o}lder inequality, the coarea formula, the isoperimetric inequality and (\ref{coarea}) to obtain
\begin{eqnarray*}4\pi |\Omega(t)|&\leq &L(t)^2=\bigg(\int_{\Gamma(t)}1\ ds_t\bigg)^2\leq\int_{\Gamma(t)}\frac{1}{|\nabla u|}\ ds_t\int_{\Gamma(t)}|\nabla u|\ ds_t\\
&=&-A'(t)\int_{\Gamma(t)}\nu^t _3\ ds_t\leq -\frac{1}{2} A'(t) |\Omega(t)| (a R^2+2b).
\end{eqnarray*}
We have utilized that  along $\Gamma(t)$, $|\nabla u|=|\nu_3^t|=\nu_3^t\geq 0$. Thus
$$8\pi\leq-A'(t)(a R^2+2b).$$
As conclusion,
$$\frac{8\pi}{a R^2+2b}\leq -A'(t),$$
and we may infer our desired estimate by an integration between $t=0$ to $t=h$ in the above inequality.
\end{proof}

\begin{corollary} Let $M$ be rotating liquid closed drop with mean curvature $2H(x)=ar^2+b$. If $a\not=0$ and $ab\geq 0$, then
the height of $M$, that is, $h=\max_{p,q\in M}|x_3(p)-x_3(q)|$, satisfies
\begin{equation}\label{area}
h\leq \frac{|a R^2+2b|}{8\pi}\ \mbox{\rm area}(M),
\end{equation}
where $R=\max_{x\in M}r(x)$.
\end{corollary}

\begin{proof}
By the Wente's symmetry result, we know that there exists a horizontal plane $P$ that is a plane of symmetry of $M$ and each one of pieces of $M$ in both sides of $P$ is a graph on $P$. For each graph, we apply the inequality (\ref{estimate}).
\end{proof}

\begin{remark} The above estimates (\ref{estimate}) and (\ref{area}) attain for cmc-surfaces. In this case, $a=0$ and $2H=b$. The estimate (\ref{estimate}) is an equality if $M$ is a spherical graph bounded by a circle. On the other hand, the inequality (\ref{area}) is an equality if $M$ is a round sphere, where $h$ is the diameter of the sphere.
\end{remark}

\section{Estimates of axisymmetric configurations}\label{sect-axi}


We consider axisymmetric stationary rotating surfaces, that is,  stationary rotating surfaces that are surfaces of revolution with respect to $x_3$-axis.
We write the generating curve $\alpha$  of a such surface as the graph of a function $u=u(r)$ and we parametrize the surface as
$x(r,\theta)=(r\cos\theta,r\sin\theta,u(r))$, $r\in[0,c)$, $\theta\in\r$. The Gauss map is
$$N(r,\theta)=\frac{1}{\sqrt{1+u'(r)^2}}(-u'(r)\cos\theta,-u'(r)\sin\theta,1).$$
With respect to this orientation,
 the mean curvature equation (\ref{laplace}) takes the form
\begin{equation}\label{eq-rotation}
\frac{u''}{(1+u'^2)^{3/2}}+\frac{1}{r}\ \frac{u'}{\sqrt{1+u'^2}}=ar^2+b.
\end{equation}
or
$$\bigg(\frac{ru'}{\sqrt{1+u'^2}}\bigg)'=r(ar^2+b).$$
A first integration   yields
$$v(r):=\frac{u'}{\sqrt{1+u'^2}}=\frac14 r(a r^2+2b)+\frac{d}{r}.$$
The integration constant $d$ describes the shape of the surface as follows. If $d\not=0$, the solutions correspond with toroidal shapes that do not intersect the $x_3$-axis. This family of surfaces have been studied, for example, in \cite{aq,gu,hc,ro}.

From now on we will restrict our analysis to the case $d=0$. Consider  initial conditions
\begin{equation}\label{laplace2}
u(0)=u_0,\hspace*{1cm}u'(0)=0.
\end{equation}
If it is necessary we indicate the dependence of the solutions with respect to the parameters as usually. Some properties of the solutions of (\ref{eq-rotation})-(\ref{laplace2}) are the following:
\begin{enumerate}
\item The existence is a consequence of standard theory.
\item A solution $u$ is symmetric with respect to $r=0$, that is, $u(-r; u_0)=u(r;u_0)$.
\item The surface is invariant by vertical displacements, that is,  if  $\lambda\in\r$, then $u(r;u_0)+\lambda=u(r;u_0+\lambda)$.
\item We have $-u(r;u_0,a,b)=u(r;-u_0,-a,-b)$. We will assume in this section that $b\geq 0$ (this does that the function $u$ is increasing near to $r=0^+$).
    \end{enumerate}
We describe the geometry of the axisymmetric closed surfaces. By the symmetry properties of the solutions of (\ref{eq-rotation})-(\ref{laplace2}), it suffices to know  the curve $u$ in the maximal interval of definition $[0,c_0)$. As we will see later, $c_0<\infty$, $u$ cannot be continued beyond $r=c_0$ but it is bounded at $r=c_0$ with derivative unbounded at the same point. Thus, the whole surface is obtained by rotating $\alpha$ with respect to $x_3$-axis and  reflecting  about the horizontal plane $x_3=u(c_0)$.

\begin{figure}[hbtp]
\includegraphics[width=6cm]{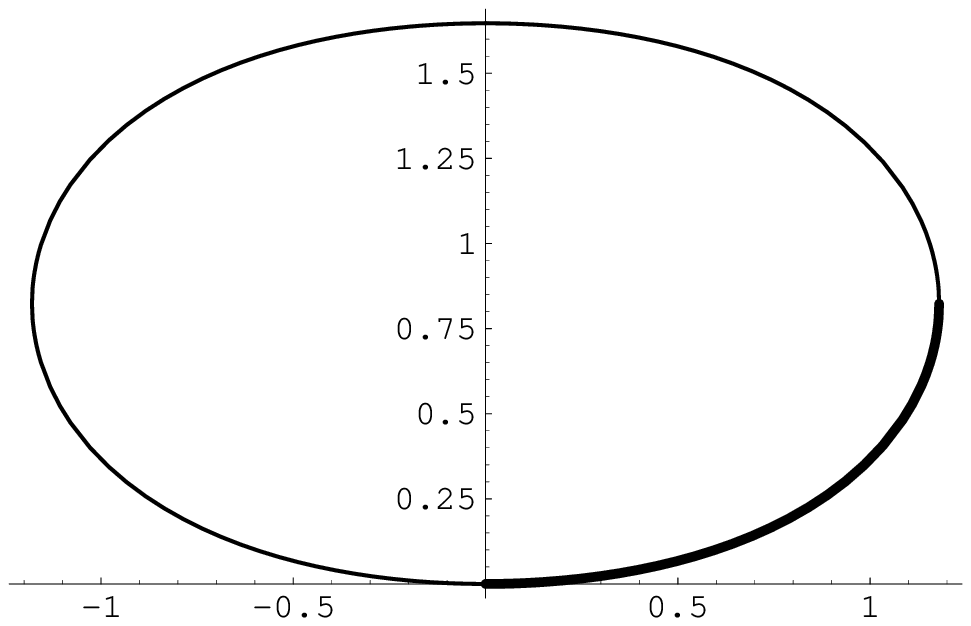}\hspace*{1cm}\includegraphics[width=6cm]{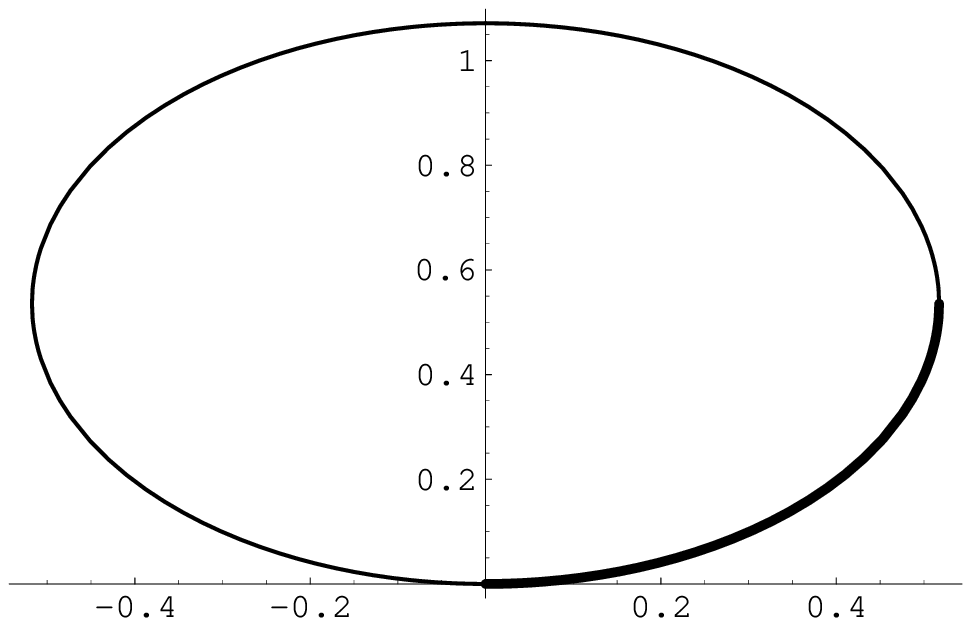}
\caption{(left) Surfaces of type I. Rotating drop for $a=1$ and $b=1$; (right) Surfaces of type II (a). Rotating drop for $a=-1$ and $b=4$.}\label{fig1}
\end{figure}

\begin{figure}[hbtp]
\includegraphics[width=6cm]{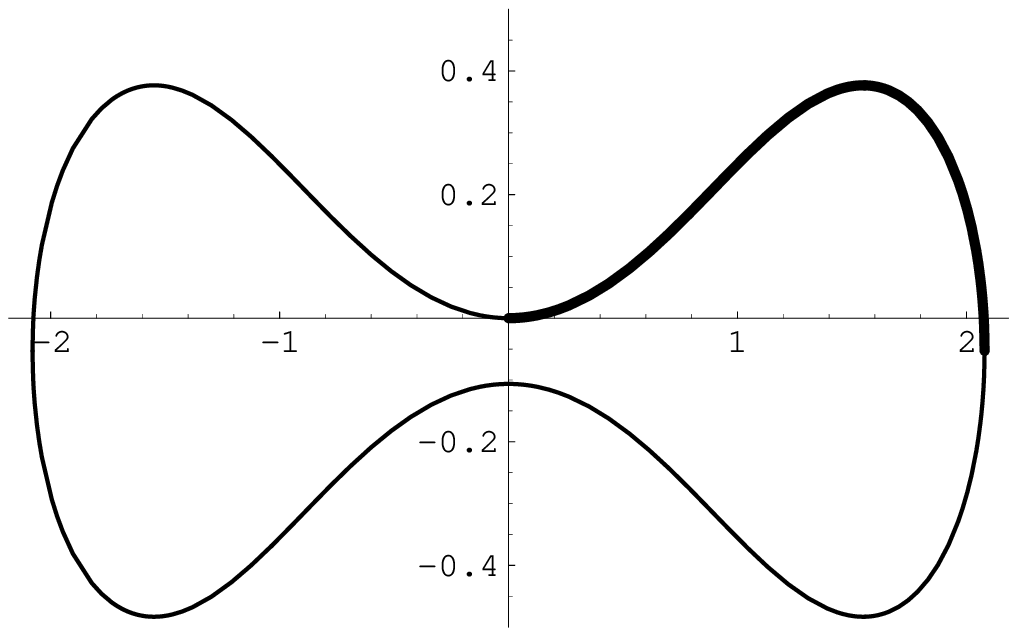}\hspace*{1cm}\includegraphics[width=6cm]{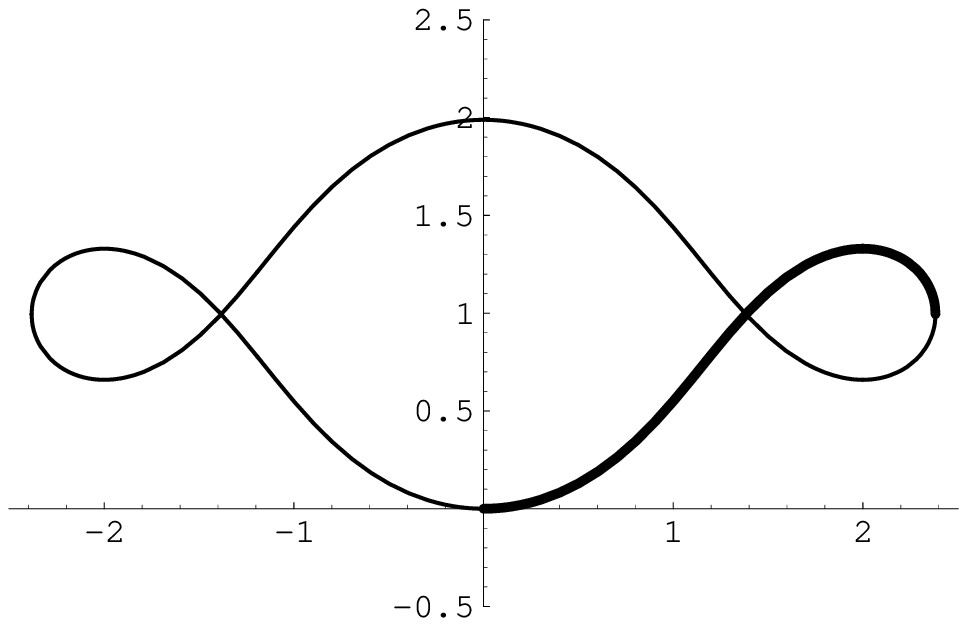}
\caption{(left) Surfaces of type II (b). Rotating drop for $a=-1$ and $b=1.2$; (right) Surfaces of type II (b). Rotating drop for $a=-1$ and $b=2$, non-embedded case.}\label{fig2}
\end{figure}

By differentiation the function $v(r)$, we obtain
\begin{equation}\label{ka}
\kappa(r):=v(r)'=\frac{u''(r)}{(1+u'(r)^2)^{3/2}}=\frac14(3a r^2+2b).
\end{equation}
Here $\kappa$ stands for the curvature of the planar curve $\alpha$. Three types of axisymmetric rotating closed drops appear and we show  the generating curves in Figures \ref{fig1} and \ref{fig2}.
The graphics correspond with  $u_0=0$ in (\ref{laplace2}). We have indicated by a bold line, the solution $u$ in the interval $[0,c_0)$. The figures  have been plotted using Mathematica (Wolfram Research Inc.).

\begin{enumerate}
\item[I]\emph{ Case $ab\geq 0$}, that is, $a>0$. Then $v(r)>0$ and $u$ is a strictly increasing function. Moreover, $\kappa(r)>0$ and $u$ is a convex function. The maximal interval where $u$ is defined is $[0,c_0)$, where $v(c_0)=1$. The surface is embedded. See Figure \ref{fig1}, left.
\item[II] \emph{ Case $ab< 0$}, that is,  $a<0$. As $\kappa(0)=b/2>0$, the function $u$ is increasing on $r$ near to $r=0^+$. The function $v(r)$ attains a maximum at $r_1=\sqrt{\frac{-2b}{3a}}$ since $v''(r_1)<0$.  Thus $u$ increasing until to reach the value  $v(c_0)=1$ if $v(r_1)\geq 1$.
     After some manipulations, this occurs iff $a\leq \frac{-2 b^3}{27}$. We distinguish  two subcases.
     \begin{enumerate}
     \item[II (a)]  Let $a\geq \frac{-2 b^3}{27}$. The function $u$ is a strictly increasing function defined in $[0,c_0)$, with $v(c_0)=1$. Moreover $u$ is a convex function. The surface is embedded. See Figure \ref{fig1}, right.
     \item[II (b)]  Let $a< \frac{-2 b^3}{27}$. The function $u$  is increasing in the interval $[0,r_2)$, with $r_2=\sqrt{\frac{-2b}{a}}$ and with an inflection at $r_1=\sqrt{\frac{-2b}{3a}}$. Next, $u$ decreases until $c_0$, where $v(c_0)=-1$. See Figure \ref{fig2}, left. The graphic of the generating curve has self-intersections iff $u_0\leq u(c_0)$. See Figure \ref{fig2}, right.
     \end{enumerate}
\end{enumerate}

We point out that for $r=c_0$, the value where $v(c_0)=\pm 1$, the function $u$ is finite. This is due to the following (\cite{hc}):
\begin{lemma} Let $f\in C^2[c_0-\epsilon,c_0]$ with $f(c_0)=\pm 1$ and $|f(r)|<1$ for $c_0-\epsilon\leq r\leq c_0$. Then the improper integral
$$\bigg|\int_{c_0-\epsilon}^{c_0}\frac{f(r)}{\sqrt{1-f(r)^2}}\ dr\bigg|$$
is finite if and only if $f'(c_0)\not=0$.
\end{lemma}
In our case, $f(r)=\frac14r(r^2+2b)$ and $f'(c_0)=\frac14(3 c_0^2+2b)\not=0$. As conclusion,
$$u(c_0)=\int_{c_0-\epsilon}^{c_0}u'(r)\ dr=\int_{c_0-\epsilon}^{c_0}\frac{f(r)}{\sqrt{1-f(r)^2}}\ dr<\infty.$$

 In this section, we compare an axisymmetric rotating drop  with  appropriate spheres. Consider the sphere obtained  by rotating about the $x_3$-axis the graphic of a function $y=y(r)$. Suppose that $u$ is solution of (\ref{eq-rotation})-(\ref{laplace2}) and that $I=[0,c]$
is an interval where $u$ is defined. In order to state our results,  we take a piece of circle
with the same slope than $u$ at $r=c$ and that coincides with $u$ at
the origin. Exactly, let
$$y(r)=R+u_0-\sqrt{R^2-r^2},\hspace*{1cm}R=c \frac{\sqrt{1+u'(c)^2}}{u'(c)}.$$
The graphic of $y$ is a piece of a lower halfcircle with $y(0)=u_0$ and $y'(0)=0$. The choice of
the radius $R$ is such that  $y'(c)=u'(c)$. Thus $y(r)$ is a solution of (\ref{eq-rotation})-(\ref{laplace2}) for $a=0$ and $b=2/R$.

\begin{theorem}\label{compara1} Let $u$ be a solution of
(\ref{eq-rotation})-(\ref{laplace2}) defined in the interval $[0,c]$. Suppose that $a>0$, $b\geq 0$. Then
$$u(r)<y(r),\hspace*{1cm}0<r\leq c.$$
\end{theorem}
\begin{proof}
Let $\psi(r)$ be the angle that makes the graphic of $u$ with the
$r$-axis at each point $r$, that is $\tan\psi(r)=u'(r)$. By the definition of the function $v(r)$,
$$\sin\psi(r)=\frac{u'(r)}{\sqrt{1+u'(r)^2}}=\frac14 r(ar^2+2b).$$
Then $\sin\psi(r)$ is positive at $(0,c]$. In particular, $u'(c)\not=0$ and the radius $R$ is well-defined. This means that $u$ is a strictly
increasing function on $r$.
The curvature $\kappa$ of $\alpha$ is also an increasing function on $r$ since
$\kappa'(r)=\frac32 a r>0$. The angle $\psi^y(r)$ and the curvature
$\kappa^y$ of the graphic of $y(r)$  are respectively
$$\sin\psi^y(r)=\frac{r}{R},\hspace*{1cm}\kappa^y(r)=\frac{1}{R}.$$
At $r=0$ we compare the curves $u$ and $y$.  We claim that
$\kappa(0)<\kappa^y(0)$. This inequality is equivalent to
\begin{equation}\label{ineb}
\frac{b}{2}<\frac{1}{R}=\frac{\sin\psi(c)}{c}=\frac{ac^2+2b}{4},
\end{equation}
which  is trivial. As $\kappa(0)<\kappa^y(0)$, $y(0)=u(0)$ and $y'(0)=u'(0)$,  the graphic of $y$
lies above of $u$ around the point  $r=0$. Theorem \ref{compara1} asserts that
this occurs in the interval $(0,c]$. Suppose, by way of contradiction, that the graphic of $u$ crosses the graphic of $y$ at some
point. Let $r=\delta\leq c$ the first value where this occurs, that is, $u(r)<y(r)$ for $r\in (0,\delta)$ and $u(\delta)=y(\delta)$.  Then $u'(\delta)\geq y'(\delta)$ and so,
$\sin\psi(\delta)\geq\sin\psi^y(\delta)$. As $u'(0)=y'(0)$, we have
\begin{equation}\label{contr1}
\int_0^\delta\bigg(\kappa(t)-\kappa^y(t)\bigg)\
dt=\int_0^\delta\bigg((\sin\psi(t))'- (\sin\psi^y(t))'\bigg)\
dt=\sin\psi(\delta)-\sin\psi^y(\delta)\geq 0.
\end{equation}
On the other hand, as  $\kappa(0)<\kappa^y(0)$ and  the above integral is non-negative, the integrand in (\ref{contr1})  is positive at some point. Then there exists $\bar{r}\in
(0,\delta)$ such that $\kappa(\bar{r})>\kappa^y(\bar{r})$.
Because $\kappa$ is increasing on $r$, we have for $r\in [\bar{r},c]$
$$\kappa(r)>\kappa(\bar{r})>\kappa^y(\bar{r})=\kappa^y(r).$$
Since $\bar{r}\leq\delta\leq c$, we have
\begin{eqnarray*}
0&<&\int_{\bar{r}}^c \bigg(\kappa(t)-\kappa^y(t)\bigg)\ dt\leq\int_{\delta}^c \bigg(\kappa(t)-\kappa^y (t)\bigg)\ dt\\
&=&\int_{\delta}^c\bigg((\sin\psi(t))'- (\sin\psi^y(t))'\bigg)\
dt=\sin\psi^y(\delta)-\sin\psi(\delta).
\end{eqnarray*}
This leads to a contradiction with (\ref{contr1}) and we have verified the theorem.
\end{proof}

For the next result, we descend vertically the circle $y(r)$ until it touches
with the graphic of $u$ at $r=c$. We call $w=w(r)$ the new position
of $y$, that is, $w(r)=y(r)-y(c)+u(c)$.

\begin{theorem}\label{compara2}
Let $u$ be a solution of (\ref{eq-rotation})-(\ref{laplace2}) defined in the
interval $[0,c]$.  Suppose that $a>0$, $b\geq 0$. Then
$$w(r)<u(r),\hspace*{1cm}0\leq r<c.$$
\end{theorem}

\begin{proof}
With a similar argument, we begin by comparying  the curvatures of $u$ and $w$ at $r=c$. Exactly, we have
$$\kappa^w(c)=\frac{1}{R}=\frac14(ac^2+2b)<\frac14(3ac^2+2b)=\kappa(c).$$
As $\kappa(c)>\kappa^w(c)$, $w(c)=u(c)$ and $w'(c)=u'(c)$, the graphic of $u$ lies above than the circle $w$ around $r=c^-$. Thus $w(r)<u(r)$ in some  interval $(\delta,c)$.
Again, the proof is by contradiction. We suppose that the
graphic of $w$ crosses the graphic of $u$ at some point. Denote by
$\delta$ the largest number such that $w(r)<u(r)$ for $r\in (\delta,c)$ and $w(\delta)=u(\delta)$.  For this value,
$w'(\delta)=y'(\delta)\leq u'(\delta)$ and $\sin\psi^y(\delta)\leq\sin\psi(\delta)$. Then
\begin{equation}\label{contr2}
\int_{\delta}^c \bigg(\kappa(t)-\kappa^w(t)\bigg) \
dt=\int_{\delta}^c\bigg((\sin\psi(t))'- (\sin\psi^y(t))'\bigg)\
dt=\sin\psi^y(\delta)-\sin\psi(\delta)\leq 0.
\end{equation}
Here we have used that $u'(c)=w'(c)=y'(c)$. As
$\kappa(c)-\kappa^w(c)>0$ and the integral in (\ref{contr2}) is non-positive, then there would be
$\bar{r}\in(\delta,c)$ such that
$\kappa(\bar{r})<\kappa^w(\bar{r})$. Because $\kappa$ is an
increasing function on $r$, for any $r\in [0,\bar{r}]$ we have
$$\kappa(r)<\kappa(\bar{r})<\kappa^w(\bar{r})=\kappa^w(r).$$
Since $\delta<\bar{r}$,
$$0>\int_0^{\delta} \bigg(\kappa(t)-\kappa^w(t)\bigg)\ dt=\int_0^{\delta}\bigg((\sin\psi(t))'- (\sin\psi^y(t))'\bigg)\
dt=\sin\psi(\delta)-\sin\psi^y(\delta),$$
where we use the fact that $u'(0)=y'(0)$.
 This contradicts  inequality (\ref{contr2}) and proves Theorem \ref{compara2}.
\end{proof}

\begin{figure}[hbtp]
\includegraphics[width=10cm]{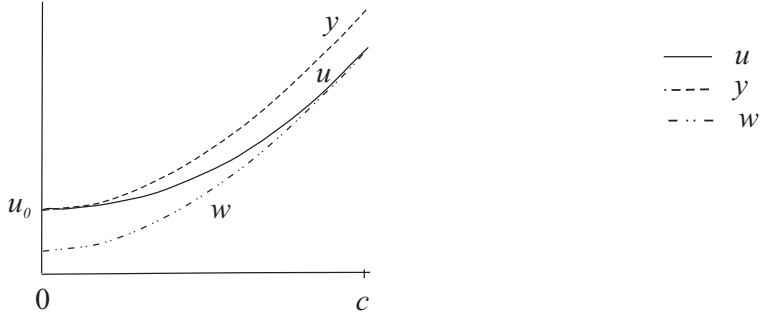}
\caption{The solution $u$ lies sandwiched  between the circles $y$ and
$w$.} \label{figcompara}
\end{figure}

As conclusion, the solution $u$ lies between two pieces of circles,
namely, $y$ and $w$, such that the slopes of the three functions
agree at the points $r=0$ and $r=c$ and the graphic of $u$ coincides with $y$ and $w$ at $r=0$ and $r=c$ respectively. See Figure \ref{figcompara}.

Finally we remark that, with appropriate modifications, the conclusions of both theorems hold even if $u$ is defined in the maximal interval $[0,c_0)$.

\begin{corollary}  Suppose that $a>0$, $b\geq 0$. Let $u$ be a solution of
(\ref{eq-rotation})-(\ref{laplace2}) defined in the interval $[0,c_0)$, where  $c_0$ is the unique positive root of $x(ax^2+2b)-4=0$.
 Then
$$u(r)-u_0<c_0-\sqrt{c_0^2-r^2},$$
\end{corollary}

An easy exercise in  calculus shows that if $a>0$, $b\geq 0$, equation $x(ax^2+2b)-4=0$ has a unique positive root.

 \begin{corollary} Let $M$ be an axisymmetric rotating liquid closed drop. Suppose that $2H(x)=a r^2+b$, where $a>0$, $b\geq 0$.  If $c_0$ is the unique  positive root of  $x(ax^2+2b)-4=0$, then the enclosed volume of the drop is less than $4/3\pi c_0^3$.
  \end{corollary}
\begin{proof} With the above notation, it is sufficient to point out that the $\mbox{volume}(u)<\mbox{volume}(w)$ and that $R=c_0$.
\end{proof}

This estimate can also be   obtained as follows. After an integration by parts, the volume of the drop is
$$V=2\pi c_0^2 u(c_0)-4\pi\int_0^{c_0}r u(r)\ dr=2\pi\int_0^{c_0}r^2 u'(r)\ dr.$$
 From the expression of $\sin\psi(r)$ and since  $ ar^2+2b\leq ac_0^2+2b$, we have
\begin{equation}\label{uprima}
u'(r)=\frac{r}{\sqrt{\bigg(\frac{4}{a r^2+2b}\bigg)^2-r^2}}<\frac{r}{\sqrt{\bigg(\frac{4}{a c_0^2+2b}\bigg)^2-r^2}}:= g(r).
\end{equation}
An explicit integration of $2\pi\int_0^{c_0}r^2 g(r)dr$ and using the fact $a c_0^3+2b c_0=4$, we conclude that  $V<\frac{4\pi}{3}c_0^3$.

 The next theorem establishes bounds for the height and the area.

\begin{theorem}\label{t-axi} Let $M$ be an axisymmetric rotating surface of type I  given by a solution $u$ of (\ref{eq-rotation})-(\ref{laplace2}) and defined in the interval $[0,c]$. Suppose that $a>0$, $b\geq 0$.  Then
\begin{equation}\label{axi1}
\frac{1}{b}\bigg(2-\sqrt{4-b^2 r^2}\bigg)\leq u(r)-u_0\leq \frac{4-\sqrt{16-r^2(a r^2+2b)^2}}{a r^2+2b},\hspace*{1cm}r\in(0,c].
\end{equation}
Denote by $A(c)$ the area of $M$. Then
\begin{equation}\label{axi2}
\frac{4\pi}{b^2}\bigg(2-\sqrt{4-b^2 c^2}\bigg)<A(c)<\frac{8\pi (4-\sqrt{16-c^2(ac^2+2b)^2})}{(a c^2+2b)^2}.
\end{equation}
\end{theorem}

\begin{proof}  The estimates are obtained by appropriate bounds for the derivative $u'(r)$. We have an upper bound for $u'$ by
 (\ref{uprima}). On the other hand, and since $ar^2+2b\geq2b$, $u'(r)\geq rb/\sqrt{4-b^2r^2}$. If we introduce both bounds of $u'(r)$ in the formulas for $u(r)$ and $A(c)$, namely,
$$u(r)-u_0=\int_0^r u'(t)\ dt,\hspace*{1cm}A(c)=2\pi\int_0^c r\sqrt{1+u'(r)^2}\ dr,$$
  the estimates follow by simple integrations. We remark that the inequality in the right hand-side of (\ref{axi1}) is also a consequence of Theorem \ref{compara1}.
\end{proof}

\begin{remark} In the case that the surface $M$ is closed, Theorem \ref{t-axi} reads as
$$u(c_0)\leq u_0+c_0,\hspace*{1cm}\mbox{area}(M)<2A(c_0)=4\pi c_0^2,$$
where $c_0$ is the unique positive root of $x(ax^2+2b)-4=0$. Here $u(c_0)-u_0$ measures the half of the distance between the highest and the lowest points of $M$. In particular, the area of $M$ is less than the area of the sphere of radius $c_0$ that contains in its inside the surface $M$: see Figure \ref{figcompara}, where the sphere determined by $w$ satisfies this property.
\end{remark}

\begin{remark} If $H$ is constant, then $a=0$ and $b=2H$ and $u$ describes a spherical cap of radius $2/|b|$. The estimates  (\ref{axi1}) are now equalities.
\end{remark}

We employ Theorem \ref{height} in the axisymmetric case.

\begin{corollary} Under the same hypothesis and notation as in Theorem \ref{t-axi}, we have
$$u(r)-u_0\leq \frac{a r^2+2b}{8\pi}A(r)< \frac{4-\sqrt{16-r^2(ar^2+2b)^2})}{a r^2+2b}.$$
\end{corollary}

\begin{proof}If the function $u$ is defined in the interval $[0,r]$, the boundary $\partial M$ of $M$ is given by the level $x_3=u(r)$. With the notation of Theorem \ref{height}, $R=r$ and $h=u(r)-u_0$, which gives the first inequality. The second one is a consequence of  (\ref{axi2}).
\end{proof}

  For cmc-graphs, a classical result due to Serrin \cite{se} asserts that if $M=\mbox{graph}(u)$ with boundary in a plane $P=\{x_3=\delta\}$, then $|u(x)-\delta|\leq 1/|H|$. This estimate can be obtained by computating the Laplacian of the function $Hx_3+N_3$. For stationary rotating graphs, we have

\begin{theorem} Let $u$ be solution of (\ref{eq-rotation})-(\ref{laplace2}) with $a,b>0$. Then
$$u(r)\leq\frac{b}{a r^2+b} u_0+\frac{2}{b}.$$
\end{theorem}

\begin{proof} Let $M$ be the corresponding  axisymmetric surface generated by $u$. From (\ref{ene2}) and (\ref{delta}), we have
\begin{eqnarray*}
\Delta(H x_3+N_3)&=&(2H^2-|\sigma|^2)N_3+2a(h-x_3 N_3)N_3\leq 2a(hx_3- N_3)N_3\\
&=&-\frac{2a t u'(t)}{1+u'(t)^2}\leq 0,\hspace*{1cm}t\geq 0,
\end{eqnarray*}
where we utilise the parametrization of $M$ as surface of revolution, that $2H^2-|\sigma|^2\leq 0$ and that $a, u'>0$. Consider $u$ defined on $[0,r]$ and thus, $\partial M$ is given by the level $u=u(r)$. Then the maximum principle yields
$$Hx_3+N_3\geq\min_{\partial M}(Hx_3+N_3)\geq \min_{\partial M}(H x_3)=H(r)u(r),$$
or, for each $t\in [0,r]$, we have
$$u(t) \geq\frac{H(r)}{H(t)}u(r)-\frac{N_3(t)}{H(t)}\geq \frac{H(r)}{H(t)}u(r)-\frac{1}{H(t)}.$$
By letting $t=0$, we obtain the desired estimate.
\end{proof}

 By numerical computations one can see that the above estimate does not hold for surfaces of type II.
We end this section showing a new use of formula (\ref{ar2}) as follows. Regarding the axisymmetric case, there exist special situations about the behaviour of the surface $M$ with respect to the plane $P$ containing $\partial M$. For example, for surfaces of type I, $M$ is  orthogonal to $P$ along $\partial M$ iff the radius of $\partial M$ is $R=c_0$ with $v(c_0)=1$.  The same occurs for surfaces of type II (a). If the surface is of type II (b), the surface is orthogonal to $P$ if $R=c_0$ with $v(c_0)=-1$. Analogously, $M$ is tangent to $P$ iff the boundary is a circle of radius $R$, with $a R^2+2b=0$ (only for surfaces of type II). We prove that this can generalize for any (non necessarily axisymmetric) stationary rotating surface.

\begin{corollary} Let $M$ be a stationary rotating surface bounded by a horizontal circle of radius $R$ centred at the $x_3$-axis. Denote by $P$ the plane containing the boundary. Let  $2H(x)=ar^2+b$ be the mean curvature of $M$.
\begin{enumerate}
\item  Suppose that $M$ is orthogonal to $P$ along $\partial M$.  Then the radius $R$ satisfies
\begin{enumerate}
\item $R(aR^2+2b)=4$ if the  boundary is running according to the counterclockwise direction.
\item $R(aR^2+2b)=-4$ if the  boundary is running according to the clockwise direction.
\end{enumerate}
\item Suppose that $M$ is tangent to $P$ along $\partial M$. Then $R=\sqrt{\frac{-2b}{a}}$. In particular $ab<0$.
\end{enumerate}
\end{corollary}
\begin{proof} It is sufficient to consider (\ref{ar2}). By distinguishing the fact that $\nu_3=\pm 1$ or $\nu_3=0$, we are going establishing the statements of the corollary.\end{proof}

\section{Rotating  liquid drops with boundary} \label{sect-embedded}

This section is devoted to  study   stationary rotating embedded surfaces with non-empty boundary.  The main tools that we will use are the maximum principle, the so-called reflection method and the integral formulae of section \ref{sect-integral}. The reflection method was employed to prove a classical result due to Alexandrov \cite{al}   that asserts that round spheres are the only embedded closed  cmc-surfaces in Euclidean space $\r^3$. The method  was used again in the  cited Wente's theorem \cite{we2}.  Our first result gives sufficient conditions to assure that the surface is a graph. As in \cite{we2}, the result holds for embedded surfaces whose mean curvature $H$ depends only on the $x_1$ and $x_2$ coordinates.

\begin{figure}[hbtp]
\includegraphics[width=13cm]{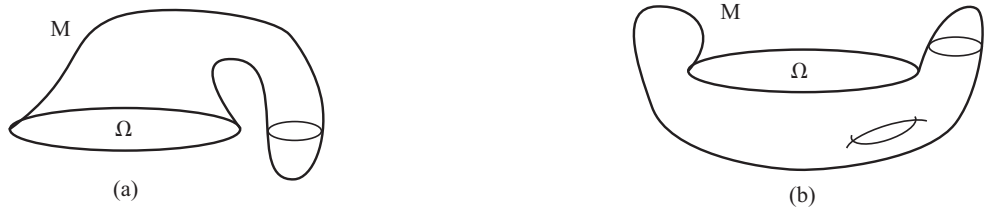}
\caption{Case (a) is forbidden by Theorem \ref{grafo}; case (b) is forbidden by Theorem \ref{t-brito}.}\label{fig33}
\end{figure}

\begin{theorem}\label{grafo}
Let $\Gamma$ be a   Jordan curve contained in a horizontal plane $P=\{x_3=d\}$ and $\Omega\subset P$ the corresponding bounded domain by $\Gamma$.
Let $M$ be an embedded  compact surface with boundary $\Gamma$ and suppose that $M$ satisfies:
 \begin{enumerate}
\item The mean curvature $H$ depends only on the $x_1$, $x_2$ coordinates.
 \item The surface $M$ does not intersect the cylinder $\Omega\times (-\infty,d]$.
 \item In a neighbourhood of $\partial M=\Gamma$ in $M$, the surface $M$ is a graph above $\Omega$.
 \end{enumerate}
 Then $M$ is a graph on $\Omega$. See Figure \ref{fig33} (a).
 \end{theorem}

\begin{proof} We apply the reflection method by using reflection with respect to horizontal planes.  For completeness, and since we will use it throughout this section, we describe the process. See \cite{we2} for details. Without loss of generality, we assume that $P$ is the plane $x_3=0$ and we define the embedded surface $T = M\cup (\Gamma\times (-\infty,0])$. This surface divides the ambient space $\r^3$ in two components. We denote by $W$ the component that contains $\Omega$ and we orient $M$ with the Gauss map $N$ that points towards $W$.

We introduce the following notation. Let $P_t$  be the 1-parameter family of translated copy of
 $P$, where we choose the parameter $t$ such that $P_t=\{x_3=t\}$. Let $A_t^+=\{x\in \r^3; x_3\geq t\}$ and $A_{t}^-=\{x\in \r^3; x_3\leq t\}$. Also, let $M_t^+=A_t^+\cap M$,
$M_t^-=A_t^-\cap M$ and $M_t^*$ the reflection of $M_t^+$ about the plane $P_t$. Because $M$ is a compact surface, for $t$ large, $P_t$ is disjoint from $M$. Now, if we approach $M$ by $P_t$ by moving down $P_t$ (letting $t\searrow 0$), one gets a  the first
plane  $P_{t_0}$, $t_0>0$, that reaches $M$, that is, $P_{t_0}\cap M\not=\emptyset$, but if $t>t_0$ then $P_t\cap M=\emptyset$.
 Thus $P_{t_0}$ is tangent to $M$ at some point and $M$ is contained in one side of $P_{t_0}$: $M\subset A_{t_0}^+$. Decreasing $t$, let consider  $M_t^*$. Since $M$ is embedded, the reflected surface $M_t^{*}$  lies inside  $W$, at least, near $t_0$: there exists at least a small $\epsilon_0>0$ such that  $M_{(t_0-\epsilon)}^{*}\subset W$ for $ 0<\epsilon<\epsilon_0$ and $M_{t_0-\epsilon}$ is a graph over $P_{t_0-\epsilon}$.

  From $t_0-\epsilon_0$ and letting $t\searrow 0$, one can  reflect $M_t$ about $P_t$, successively until one reaches a first time point of contact point of  $M_t^*$ with $M$. Exactly, consider
 $$t_1=\inf\{t<t_0;M_a^*\subset W, a\in (t,t_0]\}.$$
We claim that $t_1=0$. On the contrary, that is, $t_1>0$, $M_{t_1}^*$ and $M_{t_1}^-$ are two surfaces with $\partial M_{t_1}^*\subset \partial M_{t_1}^-=\partial M_{t_1}^*\cup\Gamma$. Moreover, $M_{t_1}^*$ and $M_{t_1}^-$ touch at an interior point $p$ or touch at a boundary point $p\in\partial M_{t_1}^*$. We remark that $p\not\in \partial\Omega$  because the surface is a graph on $\Omega$ around $\partial\Omega$ and $M_a^*\subset W$ for $t_1<a\leq t_0$. Anyway, $M_{t_1}^*$ and $M_{t_1}^-$ are one in a side of the other in a neighbourhood of $p$. As reflections invert normal vectors, the Gauss maps of both $M_{t_1}^*$ and $M_{t_1}^-$ at such point $p$ are the same. Because the mean curvature depends only the $x_1$ and $x_2$, the mean curvatures of both surfaces agree at $p$. Moreover, the mean curvature of a point of $M_{t_1}^*$ agrees with the mean curvature of the point of $M_{t_1}^-$ where vertically projects. At last, one applies either Hopf interior maximum principle or the Hopf boundary maximum principle to infer that  $M_{t_1}^*= M_{t_1}^-$ and  $P_{t_1}$ is a plane of symmetry of $M$. Because $\Gamma\subset M_{t_1}^-\setminus M_{t_1}^*$, we derive a contradiction.

As conclusion, $t_1=0$ and this means that we can go reflecting $M_t^+$ until to arrive at $t=0$, maintaining the property that
$M_t^*\subset W$ for all $t>0$. The procedure shows that in each time $t>0$, $M_t^+$ is a graph over $P_t$. This implies that $M$ is a graph on $\Omega$.
\end{proof}

In the context of this theorem,  we can obtain a similar  result as the one obtained by Wente in   \cite{we2} for the case that the surface is bounded by two curves in parallel planes. This is motivated by the physical problem of a rotating drop of liquid trapped between two parallel plates in  absence of gravity.

\begin{corollary}
Let $P_1$ and $P_2$ be two horizontal planes. Consider $\Gamma_1\cup \Gamma_2$  two Jordan curves, $\Gamma_i\subset P_i$, $i=1,2$, such that $\Gamma_2$ is the vertical translation of $\Gamma_1$ to the plane $P_2$. Denote by $\Omega$ the bounded domain determined by $\Gamma_1$ in $P_1$.
Let $M$ be an embedded  compact surface with boundary $\Gamma_1\cup \Gamma_2$ whose mean curvature depends only on the $x_1,x_2$ coordinates. Suppose that  one of the following two conditions is satisfied:
\begin{enumerate}
\item  $M$ does not intersect the solid cylinder $\Omega\times\r$.
\item $M$ is included in the solid cylinder $\Omega\times\r$.
\end{enumerate}
Then the horizontal plane $P$ equidistant from $P_1$ and $P_2$ is a plane of symmetry of $M$. Moreover,  each one of the parts of $M$ that lie in the two half-spaces determined by $P$ is a graph on $P$.
\end{corollary}

\begin{proof} Assume that $P_i=\{x_3=d_i\}$, $i=1,2$, with $d_1<d_2$. In both cases, we construct a closed surface of $\r^3$ and let us apply the reflection method with horizontal planes. It is sufficient to consider $M\cup(\Omega_0\times [d_1,d_2])$, where $\Omega_0$ is the orthogonal projection of $\Omega$ onto the plane $x_3=0$.

\end{proof}
We remark that in the above result, it is not necessary that $M$ is included in the slab determined by $P_1\cup P_2$.

Finally we use the reflection method with vertical planes in a special case.

\begin{theorem} \label{circle} Let $M$ be a rotating liquid closed drop. Assume that with the choice of the Gauss map that points inside, the mean curvature is $2H(x)=a r^2+b$, with $a<0$. Then $M$ is an axisymmetric surface (with respect to the $x_3$-axis).
\end{theorem}

\begin{proof} We prove that any plane containing the $x_3$-axis is a plane of symmetry of $M$. Without loss of generality, we are going to show that  $P_0=\{x_2=0\}$ is a plane of symmetry of $M$. For each $t\in\r$, let $P_t=\{x_2=t\}$ and consider an analogous notation as in the proof of Theorem \ref{grafo}, where $W$ is the bounded domain of $\r^3$ determined by $M$. We begin with the reflection process with planes $P_t$ and $t$ near $+\infty$ (we are assuming that $M_0^+\not=\emptyset$; on the contrary, we begin with values of $t$ near to $-\infty$). After the time $t_0>0$, we arrive the time $t=t_1$. We show that $t_1\leq 0$.  On the contrary, that is, $t_1>0$, $M_{t_1}^-$ and $M_{t_1}^*$ have a common (interior or boundary) contact point $p$, where $M_{t_1}^*$ locally lies over $M_{t_1}^-$ with respect to the vector $N(p)$. Assume that $p=q^*$ with $q\in M_{t_1}^+$, the reflection of the point $q$ about the plane $P_{t_1}$. We have two cases:
\begin{enumerate}
\item If $p$ is an interior point of both surfaces $M_{t_1}^-$ and $M_{t_1}^*$, then $H(q)\geq H(p)$, that is, $ar(q)^2\geq a r(p)^2$. As $a<0$, $r(q)\leq r(p)$: contradiction (even if $p\in M_0^-$).
\item If $p\in\partial M_{t_1}^-\cap\partial M_{t_1}^*$, then $q^*=q=p$. The maximum principle says that $H(x)\geq H(y)$ for $x\in M_{t_1}^*$, $y\in M_{t_1}^-$ near to $p$ with $\pi(x)=\pi(y)$, being $\pi$ the orthogonal projection onto the common tangent plane $T_p M_{t_1}^-=T_p M_{t_1}^*$. Then  $a r(z)^2\geq a r(y)^2$, with $x=z^*$, $z\in M_{t_1}^+$, that is, $r(z)\leq r(y)$, obtaining  a contradiction again.
\end{enumerate}
As conclusion, we have proved that $t_1\leq 0$. In particular, $M_{t_1}^-\not=\emptyset$. With a similar reasoning and  with vertical planes coming from $t=-\infty$, we show that $t_1\geq 0$ and thus $t_1=0$ proving that $P_0$ is a plane of symmetry of $M$.
\end{proof}

With a s similar argument, we obtain the version of Theorem \ref{circle} in the case that $\partial M\not=\emptyset$.

\begin{corollary} Let $P$ a horizontal plane and let $Q$ be a plane containing the $x_3$-axis. Set $R=P\cap Q$. Consider $\Gamma\subset P$ a Jordan curve and denote by $\Omega\subset P$ the bounded domain by $\Gamma$. Assume that $\Gamma$ is symmetric with respect to the reflections about $Q$ and that $R$ divides $\Gamma$ into two pieces that are graphs on $R$. Let  $M$ be a rotating liquid drop bounded by $\Gamma$ such that
$M$ lies in one side of $P$ and that with the choice of the Gauss map that points inside of $M\cup\Omega$, the mean curvature $2H(x)=ar^2+b$ satisfies $a<0$. Then $Q$ is a plane of symmetry of $M$. In the particular case that $\Gamma\subset P$ is a circle centred at that $x_3$-axis, then  $M$ is an axisymmetric surface.
\end{corollary}

The following result compares rotating graphs with  stationary rotating surfaces that lie in solid vertical cylinders.

\begin{theorem} Let $\Gamma$ be a   Jordan curve contained in a horizontal plane $P$ and $\Omega$ the corresponding bounded planar domain.  Suppose that there exists a stationary rotating graph $G$ on $\Omega$ with $\partial G=\Gamma$ and assume that the mean curvature $2H_G(x)=ar^2+b$ does not vanish at any point of $G$.  Let $M$ be a stationary rotating surface bounded by  $\Gamma$ with the following conditions:
\begin{enumerate}
\item The mean curvature $H$ of $M$ satisfies  $|2H(x)|=|ar^2+b|$.
\item The surface  $x(M)$ lies in  $\Omega\times\r$.
\end{enumerate}
Then $x(M)$ either coincides with $G$ or with its reflection about $P$.
\end{theorem}

\begin{proof} Denote by $G^*$ the reflection of $G$ about the plane $P$.  Without loss of generality, we suppose that the mean curvature $H_G$ is positive with the orientation $N_G$ on the graph $G$ pointing downwards (by the maximum principle, this implies that $G$ lies above the plane $P$). We choose  the orientation $N$ on $M$ whose mean curvature is exactly $H_G$.  We move $G$ upwards so it does not touch $M$ and then we drop it until it reaches a contact point $p$ with $M$ for the first time. Denote by $G'$ the translated graph at this time. If $p\not\in \Gamma$, then $M$ and $G'$ are tangent at $p$ and  $N(p)=\pm N_G(p)$. If $N(p)=N_G(p)$, and as   the mean curvatures of $G'$ and $M$ agree for the same choice of normal vector fields at $p$, the maximum principle implies that $M$ and $G'$ should be coincide. This is a contradiction, since the boundaries $\partial M=\Gamma$ and $\partial G'$ lie at different heights. If $N(p)=-N_G(p)$, we change the orientation of $M$, namely, $N'=-N$, and the mean curvature  $H'=-H<0$. Then $N'(p)=N_G(p)$, but $H'(p)<H_G(p)$:  by using the maximum principle, we arrive to a contradiction.

As conclusion, we can move $G$ downwards until that $G$ returns into its original position. Working now with the graph $G^*$, the same reasoning shows that $M$ lies above $G^*$. If $\nu_M$ and $\nu_G$ denote the inner conormal unit vectors of $M$ and $G$ respectively along their common boundary $\Gamma$, we have proved that
$$|\langle \nu_M(p),E_3\rangle|\leq\langle\nu_G(p),E_3\rangle,\hspace*{1cm}\forall p\in\Gamma.$$
If the equality holds at some point $p_0\in\Gamma$, this means that $M$ is tangent to $G$ or $G^*$ at $p_0$. Now we use the boundary maximum principle which implies  that either $M=G$ or $M=G^*$, proving the result.

Suppose, by way of contradiction, that we have
$$|\langle \nu_M(p),E_3\rangle| <\langle\nu_G(p),E_3\rangle$$
for each $p\in\Gamma$. Integrating this  inequality along the common boundary $x(\partial M)=\Gamma=\partial G$, we have
\begin{equation}\label{cylinder}
\bigg|\int_{\partial M}\langle \nu_M,E_3\rangle\ ds\bigg|\leq\int_{\partial M}|\langle \nu_M,E_3\rangle|\ ds <\int_{\partial G}\langle \nu_G,E_3\rangle\ ds.
\end{equation}
We employ the integral formula (\ref{ar2}) for both surfaces $M$ and $G$:
$$\int_{\partial M}\langle \nu_M,E_3\rangle\ ds=-\frac14\int_{\partial M}(ar^2+2b)\langle \alpha\wedge \alpha',E_3\rangle\ ds.$$
$$\int_{\partial G}\langle \nu_G,E_3\rangle\ ds=-\frac14 \int_{\partial G}(ar^2+2b)\langle \alpha\wedge \alpha',E_3\rangle\ ds.$$
Since the integrands in the right-hand sides of the above two equations depend only on $\Gamma$ and they are the same or of reverse sign depending of the direction of $\alpha_{|\partial M}$ and $\alpha_{|\partial G}$, we derive a contradiction with the inequality (\ref{cylinder}).
\end{proof}

We now extend Theorem \ref{height} for rotating liquid drops with boundary.

\begin{theorem} Let $P$ be a horizontal plane of $\r^3$. Consider $M$  a stationary rotating embedded surface with $\partial M\subset P$. Assume that the mean curvature is $2H(x)=ar^2+b$, where $a\not=0$ and $ab\geq 0$.  If $M$ lies in one side of $P$, then
$$\max_{x\in M}\mbox{\rm dist}(x,P)\leq\frac{|a R^2+2b|}{8\pi}\ \mbox{ \rm area}(M),\hspace*{1cm}R=\max_{x\in M}r(x).$$
\end{theorem}

\begin{proof} Without loss of generality, we suppose that $M$ lies in the upper half-space determined by $P$. Consider $W$ the domain that encloses $M\cup\Omega$, being $\Omega$  the bounded planar domain by $\partial M$. We orient $M$ by the Gauss map that points inside $W$. By considering the highest point of $M$ and the maximum principle, we deduce that the mean curvature is positive. Thus $a,b>0$. We use the reflection method with horizontal planes. With the notation of Theorem \ref{grafo}, we have the next possibilities at the time $t_1$:
\begin{enumerate}
\item  $t_1=0$. Then $M$ is a graph on $\Omega$ and  Theorem \ref{height}  proves the result.

\item $t_1>0$ and there exists a tangent point between the surfaces $M_{t_1}^{-}$ and $M_{t_1}^*$. Then the maximum principle says that $P_{t_1}$ is a plane of symmetry of $M$, which it is a contradiction because the boundary $\partial M$ lies below $P_{t_1}$.
\item $t_1>0$ and  $M_{t_1}^*$ and $\partial M$ contact at some point. Then $M_{t_1}^{+}$ is a graph over a domain of $P_{t_1}$. Theorem \ref{height} says again
\begin{eqnarray*}
\frac12 \max_{x\in M}\mbox{\rm dist}(x,P)&\leq& \max_{x\in M_{t_1}^+}\mbox{\rm dist}(x,P_{t_1})\leq\frac{|a R^2+2b|}{8\pi}\ \mbox{ \rm area}(M_{t_1}^+)\\
&\leq& \frac{|a R^2+2b|}{8\pi}\ \mbox{ \rm area}(M),
\end{eqnarray*}
and the desired estimate is obtained again.
\end{enumerate}
\end{proof}

We end this section with a result motivated by what happens in the theory of cmc-surfaces with boundary. If $M$ is a  cmc embedded surface with boundary, one asks  under what conditions the symmetries of the boundary of a cmc-surface are inherited by the whole surface. For example, if $\Gamma$ is a circle and $\partial M=\Gamma$, is $M$ a spherical cap?   If $M$ lies in one side of the plane containing $\Gamma$,  the Alexandrov reflection method proves that $M$ is a spherical cap. Thus, one seeks conditions that assure that the surface lies in one side of $P$. Two results stand out  in this setting. Assume that $\Gamma$ is a closed curve contained in a plane $P$. The first one is due to Koiso \cite{koi} and shows that  if  $M$ does not intersect the outside of $\Gamma$ in $P$, then $M$ lies in one side of $P$. The second result, due to Brito, Sa Earp, Meeks and Rosenberg \cite{bemr}, shows that if $\Gamma$  is strictly convex
and $M$ is  transverse to $P$  along the boundary $\partial M$, then $M$ is entirely
contained in one of the half-spaces of $\r^3$ determined by $P$.  Here, transversality means that the surface $M$ is never tangent to the
plane $P$ along $\partial M$.

\begin{figure}[hbtp]
\includegraphics[width=12cm]{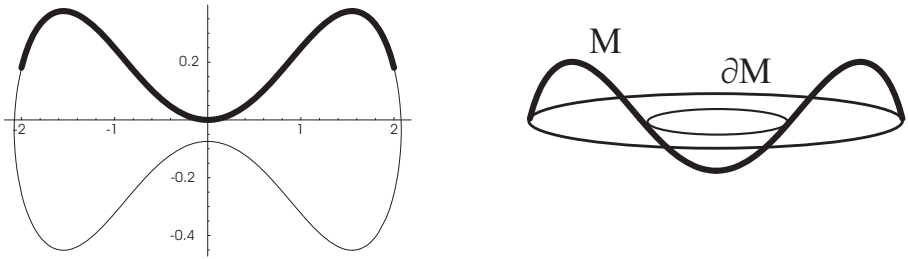}
\caption{Rotating drop for $a=-1$ and $b=1.2$. The bold line represents the generating curve of a rotating drop that is embedded and  intersects the outside of the boundary in the plane $P$ containing the boundary. Moreover, the surface is transverse to $P$. Anyway, the surface does not lie in one side of $P$.}\label{fig3}
\end{figure}

For stationary rotating surfaces, one poses the same question whether a stationary rotating surface bounded by a (horizontal) circle  must be a surface of revolution. A first difference is that even if $M$ lies in one side of $P$, $P$ the plane containing $\partial M$, one cannot apply the Alexandrov reflection method  (a special case is Theorem \ref{circle}).  With respect to Koiso's theorem, actually her result holds assuming that the mean curvature does not change of sign and thus, it is true for stationary rotating surfaces where $H$ does not vanish.   Finally, the result cited in \cite{bemr} does not hold for stationary rotating surfaces even if the surface is axisymmetric:  surfaces of type II (b) provide a counterexample of both results for cmc-surfaces as it can see in Figure \ref{fig3}. The last result in this section is related with this theorem.

\begin{theorem} \label{t-brito}Let $\Gamma$ be a closed curve  contained in a horizontal plane $P$ and star-shaped with respect to $O$, the intersection point between $P$ and the $x_3$-axis. Let $\Omega\subset P$ be the corresponding bounded planar domain by $\Gamma$.  Let $M$ be a stationary rotating  embedded   surface  with boundary $\Gamma$ and suppose that $M$ satisfies:
 \begin{enumerate}
\item The mean curvature does not vanish on $M$.
 \item The surface $M$ does not intersects the domain $\Omega$.
 \item The surface $M$ is transverse to $P$ along $\Gamma$.
 \end{enumerate}
 Then $M$ lies in one side of $P$. In particular, the result holds for horizontal convex curves whose inside intersects the $x_3$-axis.
 See Figure \ref{fig33} (b).
\end{theorem}

\begin{proof}
Without loss of generality, we suppose that $P$ is the plane $x_3=0$ and $H$ is positive on $M$. By transversality,  in a neighborhood of the boundary $\partial M$ the surface $M$ is contained in one of the two connected components of $\r^3\setminus P$, which, without loss of generality, can be assumed to be, the upper half-space. In this situation, we will prove that $M$ is above $P$.  We attach the domain $\Omega$ to $M$, obtaining a closed surface $M'=M\cup\Omega$. Thus, $M'$ encloses a domain $W$ of $\r^3$. We orient $M'$ by the mean curvature vector of $M$ and let us denote by $\eta_{\Omega}$ the induced orientation on $\Omega$.

We claim that $N'$ points to the domain $W$. For this, we take the highest point $p$ of $M$ with respect to the plane $P$. In particular, $p$ is an interior point of $M$ and $N(p)=\pm E_3$. As the mean curvature $H$ of $M$  is positive, the maximum principle implies that $N(p)=-E_3$, and so, $N$ points to $W$.

We show the theorem by way of contradiction: assume that $M$ has points below $P$. We  use Equation (\ref{ar2}) and the notation that appears there. Because $\Gamma $ is a star-shaped curve with respect to the origin, the function $\langle x\wedge\alpha',E_3\rangle$ has sign along $\partial\Omega$. We show  that this sign is positive.  Since $M$ does not intersect $\Omega$ and $N$ points towards $W$,   the orthogonal projection of the restriction $N_{|\partial M}$  onto the plane $P$  points outside $\Omega$. This means that $\alpha'=\nu\wedge N_{|\partial\Omega}$ follows the counterclockwise direction along $\partial M$. Thus $\langle x\wedge\alpha',E_3\rangle>0$ along $\partial M$. Once proved this, we remark that the function $\nu_3=\langle \nu,E_3\rangle$ is also positive.
Using Equation (\ref{ar2}), we arrives to a contradiction and this completes the proof.
\end{proof}

\begin{remark} Actually, one can replace the hypothesis on the transversality by the fact that  the surface lies locally in one side of $P$ around $\Gamma$. With the same notation as in the above proof, this means that  $\int_{\partial M}\nu_3\ ds\geq 0$. Thus, the left side of (\ref{ar2}) implies that $a=b=0$ and $H=0$: contradiction.
\end{remark}


\section{Stability}\label{sect-estable}


In this section we give two results on stability. The first one refers to strongly stability and assures that any stationary rotating graph is strongly stable. This generalizes a well known result of the theory of minimal surfaces in Euclidean space. According to Section \ref{sect-preli}, we define the energy of a stationary rotating surface as follows. Let $x:M\rightarrow \r^3$ be an immersed compact surface with mean curvature $2H(x)=a r^2+b$ for any $x\in M$. We define the energy of the immersion $x$ as
$$E(x)=\int_M 1\ dM+a\int_M  r^2 x_3 N_3\ dM+ b\int_M x_3 N_3\ dM.$$
The description of each one of the integrals that appears in the right-hand side is the following. The first one represents the area of the surface and is proportional to the surface tension energy; the second integral is the energy of the centrifugal force of the surface with respect to the $x_3$-axis;  and the last one is the algebraic volume between the surface and the plane $x_3=0$.

\begin{theorem}
Let $M$ be  a stationary rotating surface that is a graph on some horizontal domain $\Omega\subset\r^2$. Then $M$ is strongly stable. Moreover, there  holds  the following property about the energy of $M$. Let $y:M'\rightarrow\r^3$ be an immersion of a compact oriented surface $M'$ with the same boundary as $M$ in such way that $M\cup M'$ defines an oriented 3-chain $W$. If $M'$ is included in the vertical solid cylinder $\Omega\times\r$, then $E(M)\leq E(y)$.
\end{theorem}

\begin{proof}

The function $N_3=\langle N,E_3\rangle$ is a non-vanishing function on the surface that satisfies $L[N_3]=0$: see (\ref{ene2}). It follows from standard theory that this is equivalent to say that $M$ is strongly stable.

For the proof of the second part of Theorem, we use an argument  of calibration type . Assume that $M=\mbox{graph}(u)$, where $u$ is a smooth function defined on $\Omega$. Consider the orientation on $M$ pointing upwards, that is,
$$N(x)=\frac{1}{\sqrt{1+|Du|^2}}\bigg(-\frac{\partial u}{\partial x_1},-\frac{\partial u}{\partial x_2},1\bigg)(x),\hspace*{1cm}x\in M$$
 On $\Omega\times\r$, we define the vector field
$$Z(x_1,x_2,x_3)=N(x_1,x_2)+(a r^2+b  )x_3E_3,$$
where the mean curvature of $M$ is $H=ar^2+b$. Then
$$\mbox{div}_{\r^3}Z=\mbox{div}_{\r^3}N+(ar^2+b)=-2H+(ar^2+b)=0.$$
Using the divergence theorem in the chain $W$, one has
$$0=\int_M\langle Z,N\rangle\ dM-\int_{M'}\langle Z,N'\rangle\ dM',$$
where $N'$ is the orientation on $M'$ induced by $W$. The first integral in the right-hand side is
$$\int_M\langle Z,N\rangle\ dM= \int_M\bigg(1+(a r^2 +b )x_3N_3\bigg)\ dM=E(M).$$
On the other hand, the second integral is
\begin{eqnarray*}
\int_{M'}\langle Z,N'\rangle\ dM' &=& \int_{M'}\bigg(\langle N,N'\rangle+(a r^2 +b )x_3N'_3\bigg)\ dM'\\
&\leq & \int_{M'}\bigg(1+(a r^2 +b )x_3N'_3\bigg)\ dM'=E(y).
\end{eqnarray*}
where here we have used $\langle N,N'\rangle\leq 1$. Hence it follows the result.
\end{proof}

The second result is about the stability of axisymmetric rotational closed surfaces. Let $M$ be a such surface whose mean curvature is $2H(x)=ar^2+b$ and let $N$ be its Gauss map.    From (\ref{n1}) and because $\partial M=\emptyset$, the coordinates functions of $N$ satisfy $\int_M N_i\ dM=0$. Thus they are test functions to study the stability of the surface. Consider $i\in\{1,2\}$. Equation (\ref{ene}) leads to $L[N_i]=-2a x_i$. We then compute  $-\int_M N_i\cdot L[N_i]\ dM$:
$$ -\int_M N_i \cdot L[N_i]\ dM=2a\int_M x_i N_i\ dM,\hspace*{1cm} i=1,2.$$
We calculate this integral on $M$.
Assume that $M$ is symmetric with respect to the plane $x_3=0$. We parametrize the lower part of $M$, that is, $M\cap \{x_3\leq 0\}$, as
$x(r,\theta)=(r\cos\theta,r\sin\theta,u(r))$, $r\in[0,c_0)$, $\theta\in\r$. Here $u=u(r)$ is a solution of (\ref{eq-rotation})-(\ref{laplace2}) and $[0,c_0)$ is the maximal interval of definition. For our purposes, it is sufficient to consider the function $N_1$. Then
$$2a\int_M x_1N_1\ dM=-4a \int_0^{c_0} \int_0^{2\pi} r^2\cos^2\theta u'(r)\ dr\ d\theta=-4\pi a\int_0^{c_0} r^2 u'(r)\ dr.$$
The surfaces of type I satisfy $a>0$ and  $u'(r)>0$. Thus

\begin{theorem} Axisymmetric rotating surfaces of type I are not  stable.
\end{theorem}

We end this paper with  several natural  questions that could and should be addressed within this theory of stationary rotating surfaces.

\begin{enumerate}

\item  Let $M$ be a rotating liquid drop with non-empty boundary. Assume that $\partial M$ lies in a horizontal plane and that the mass center of $\partial M$ lies in the $x_3$-axis. Does the mass center of $M$ lies in the $x_3$-axis? For example, if $\partial M$ is a horizontal circle centred at $x_3$-axis.
  \item Let $M$ be a rotating liquid drop whose boundary is a circle in a horizontal plane. Is $M$ a surface of revolution? The same if the boundary are two coaxial circles in horizontal planes.
  \item What  axisymmetric stationary rotating surfaces bounded by a circle are stable?
\end{enumerate}

\end{document}